\documentclass[12pt]{amsart}

\usepackage{amssymb, amsthm, amsmath, amsrefs}
\usepackage{fullpage}

\newtheorem{theorem}{Theorem}

\newtheorem{corollary}{Corollary}
\newtheorem{question}{Question}
\newtheorem*{thma}{Theorem A}
\newtheorem*{thmastar}{Theorem A*}
\newtheorem*{thmc}{Theorem C}
\newtheorem*{thmb}{Theorem B}

\theoremstyle{definition}

\theoremstyle{remark}

\numberwithin{equation}{section}

\newcommand{\D}{\mathbb{D}}
\newcommand{\R}{\mathbb{R}}
\newcommand{\cD}{\overline{\D}}
\newcommand{\E}{\mathbb{E}}
\newcommand{\C}{\mathbb{C}}
\newcommand{\T}{\mathbb{T}}

\newcommand{\rank}{{\rm rank }\,}
\newcommand{\mcP}{\mathcal{P}}

\title[Determinantal representations]{Determinantal representations of
  semi-hyperbolic polynomials}

\author{Greg Knese}

\address{Washington University in St. Louis, Department of
  Mathematics, St. Louis, MO, 63130}

\date{\today}

\email{geknese@wustl.edu}


\thanks{This research was supported by NSF grant DMS-1363239}


\begin{document}
\bibliographystyle{plain}
\maketitle

\begin{abstract} 
  We prove a generalization of the Hermitian version of the
  Helton-Vinnikov determinantal representation for hyperbolic
  polynomials to the class of \emph{semi-hyperbolic} polynomials, a
  strictly larger class, as shown by an example. We also prove that
  certain hyperbolic polynomials affine in two out of four variables
  divide a determinantal polynomial. The proofs are based on work
  related to polynomials with no zeros on the bidisk and tridisk.
\end{abstract} 

\section{Introduction}

A homogeneous polynomial $P \in \R[x_0,x_1,\dots, x_{n}]$ is
\emph{hyperbolic} of degree $d$ with respect to $e \in \R^{n+1}$ if
$P(e) \ne 0$ and if for all $x\in \R^{n+1}$ the one variable
polynomial $t \mapsto P(x-te)$ has only real zeros.  This concept was
originally studied by G\r{a}rding for its relation to PDE (see
\cite{Garding}, \cite{HL}) but it---and the related concept of stable
polynomials---has since become important to convex optimization,
combinatorics, probability, combinatorics, and analysis.  See the
papers and surveys \cite{Renegar}, \cite{Gurvits}, \cite{Wagner},
\cite{Pemantle}, \cite{HL}, \cite{MSS}. 

A deep result in the area is a determinantal representation for
trivariate hyperbolic polynomials due to Helton-Vinnikov \cite{HV}, \cite{vV}
which solved a 1958 conjecture of Lax \cite{pL} (see \cite{LPR}) and,
as is mentioned in \cite{HL}, can be used to develop the full
G\r{a}rding theory of hyperbolicity.  

\begin{thma} Let $p \in \R[x_0,x_1,x_2]$ be hyperbolic of degree $d$
  with respect to $e_2$ and monic in $x_2$.  Then, there exist
  $d\times d$ real symmetric matrices $A_0,A_1$ such that
\[
p(x_0,x_1,x_2) = \det (x_0A_0+x_1A_1+x_2I).
\]
\end{thma}

If we relax the problem to finding self-adjoint matrices instead of
real symmetric matrices, proofs more amenable to computations are
possible (see \cite{GKVW},\cite{PV},\cite{vV12}).  The resulting
theorem is just as useful for most purposes.

\begin{thmastar}
  Let $p \in \R[x_0,x_1,x_2]$ be hyperbolic of degree $d$ with respect
  to $e_2$ and monic in $x_2$.  Then, there exist $d\times d$
  self-adjoint matrices $A_0,A_1$ such that
\[
p(x_0,x_1,x_2) = \det (x_0A_0+x_1A_1+x_2I).
\]
\end{thmastar}

Our immediate goal is to prove a generalization of this result based
on a result in Geronimo et al \cite{GIK} and an extension to four
variables based on a result in Bickel and Knese \cite{BK}, while our
larger goal is to advertise the close connection between determinantal
representations of hyperbolic polynomials and sums of squares
decompositions for multivariable Schur stable polynomials.  See
\cite{GKVWdetreps}, \cite{GKVWscatter}, \cite{gKsacrifott},
\cite{gKrifitsac} for background on the latter topic.

Our main result establishes a determinantal representation with the
assumption of hyperbolicity weakened.  We shall call a homogeneous
polynomial $P \in \R[x_0,x_1,\dots, x_{n}]$ a \emph{semi-hyperbolic
  polynomial} with respect to the direction $e \in \R^{n+1}\setminus
\{0\}$ if for every $x \in \R^{n+1}$ the univariate polynomial $t
\mapsto P(x-te)$ is either identically zero or only has real roots.
The key distinction between hyperbolic and semi-hyperbolic polynomials
is that we do not assume $P(e)\ne 0$.  Some references actually
confuse the two, while Renegar \cite{Renegar} is the only reference we
have found that emphasizes the distinction.  We elaborate on our
motivations in Section \ref{conclusion}.  We do need to allow for
$t \mapsto P(x-te)$ to be identically zero, because for instance if
$P(e) = 0$ and $x=0$, then $P(-te) \equiv 0$.  We give an example of a
semi-hyperbolic polynomial that is not hyperbolic in any direction in
Section \ref{sec:ex}.  

Here is our main theorem.

\begin{theorem} \label{mainthm} Let $p \in \R[x_0,x_1,x_2]$ of degree
  $d$ be semi-hyperbolic with respect to $e_2=(0,0,1).$ Then, there
  exist $d\times d$ self-adjoint matrices $A_0,A_1,A_2$ with $A_2$
  positive semi-definite and a constant $c \in \R$ such that
\[
p(x) = c\det(x_0 A_0+x_1A_1 + x_2 A_2).
\]
Assuming $p$ has no factors depending on $x_0,x_1$ alone, the above
data can be chosen to additionally satisfy:
\begin{itemize}
\item
\[
\rank A_1 = \deg_1 p, \quad \rank A_2 = \deg_2 p,
\]
\item 
$A_1 = B_{+}-B_{-}$ with $B_{\pm}$ both positive semi-definite where
 $\rank B_{-}$ equals the number of roots of $p(1,t,i)$ in the upper half
plane and $\rank B_{+} + \rank B_{-} = \rank A_1$,
\item and $B_{-} + B_{+} + A_2 = I$.
\end{itemize}
\end{theorem}

See Section \ref{sec:btomain} for the proof.
We can recover Theorem A* when $p(e_2)\ne 0$ since $p$ will then have
degree $d$ in $x_2$ and then $A_2$ will be positive definite.  We can
then factor $A_2^{1/2}$ from the right and left of $\sum_{j=0}^{2}
x_jA_j$ in order to get a determinantal representation of the form
given in Theorem A*, namely with $A_2=I$.

There is nothing special about the vector $e_2$; a linear change of
variables could be used to establish a determinantal representation
for other semi-hyperbolic polynomials.  The assumption of no factors
depending on only $x_0,x_1$ is there to avoid certain annoyances that
such trivial factors introduce.  For instance $p(x_0,x_1,x_2) = x_1$
is certainly semi-hyperbolic in the direction $e_2$ but then
$A_2=A_0 = 0$ and the signature of the $1\times 1$ matrix $A_1$ does
not really provide any useful information.

It follows that a trivariate semi-hyperbolic polynomial
$p$ can be lifted to a four variable polynomial
\[
P(x_0,x_1,y_1,x_2) = c \det(x_0A_0 + x_1B_{+} + y_1 B_{-} +x_2 A_2)
\]
which is hyperbolic in the direction $(0,1,1,1)$ and
$P(x_0,x_1,-x_1,x_2) = p(x_0,x_1,x_2)$.  So, we are projecting a
hyperbolic polynomial (possessing a definite determinantal
representation) of four variables to a set where it is not necessarily
hyperbolic.  It also follows that a trivariate semi-hyperbolic
polynomial is a limit of hyperbolic polynomials.  Indeed, writing $p$
as in Theorem \ref{mainthm}, define for $\epsilon >0$
\begin{equation} \label{pepsilon}
p_{\epsilon}(z) = c \det(x_0 A_0+x_1 A_1 + x_2(A_2+\epsilon I)).
\end{equation}
Then, $p_\epsilon \to p$ as $\epsilon \searrow 0$.  We do not know if
semi-hyperbolic polynomials in more than three variables are the limit
of hyperbolic polynomials.

The theorem above has an curious asymmetry in its treatment of $x_0$
and $x_1$.  This is partly due to idiosyncrasies of our proof but we
also think there are some subtleties to resolve.  To be specific, one
could certainly break up $A_0$ into a difference of positive
semi-definite matrices according to its signature but we have been
unable to connect the signature of the $A_0$ we construct with
geometric properties of $p$.  We have no reason to believe this cannot
be done, especially because this issue does not arise in the
hyperbolic case.  Indeed, we can take $A_2=I$ and the number of zeros
of $t\mapsto p(tx_0,t x_1, i)$ in $\C_{+}$ equals the number of
negative eigenvalues of $x_0A_0+x_1A_1$.  Similarly, the number of
zeros of $t\mapsto p(tx_0,t x_1,-i)$ equals the number of positive
eigenvalues of $x_0 A_0+x_1A_1$.  Thus, the signature of
$x_0A_0+x_1A_1$ can be derived from properties of $p$ in the
hyperbolic case.  Notice that we evaluate $p$ on the complex line
$(0,t, i)$ to determine the signature of $A_1$ while in the theorem
above we evaluate on the line $(1,t,i)$, which actually seems less
natural.  The example in Section \ref{sec:ex} shows this is actually
necessary: using the line $(1,t,i)$ we get a correct count of the
negative eigenvalues of $A_1$ while using the line $(0,t,i)$ we get an
incorrect count.  The details are recorded in Section \ref{sec:ex}.
 
As a nice corollary, we can quickly recover the following variant of
Theorem A*.  The original proof, while not difficult, requires
transforming a real stable polynomial to a hyperbolic polynomial
through a linear transformation.  Our signature count of $A_1$ in
Theorem \ref{mainthm} makes the proof go smoothly.

\begin{corollary}[See Theorem 6.6 of \cite{BB}] \label{cor}
  If $p \in \R[x_0,x_1,x_2]$ is homogeneous of degree $d$ and
  hyperbolic with respect to all vectors in the cone $\{(0,v_1,v_2):
  v_1,v_2 >0\}$, then there exist $d\times d$ self-adjoint matrices
  $A_0,A_1,A_2$ and a constant $c\in \R$ such that $A_1,A_2$ are
  positive semi-definite, $A_1+A_2=I$, and 
\[
p(x) = c \det (x_0 A_0 + x_1A_1+x_2 A_2).
\]
\end{corollary}

Since \cite{BB} uses Theorem A to prove the above result, all of the
matrices can be taken to be real but our proof does not yield this.
For $p$ as in the corollary, $p(1,x_1,x_2)$ is known as a \emph{real
  stable} polynomial.  This formula was used in the recent paper
regarding the Kadison-Singer problem \cite{MSS}.  See Section
\ref{sec:cor} for the very short proof of the corollary.

The key tool for the proof of Theorem \ref{mainthm} is a determinantal
representation proven in Geronimo-Iliev-Knese \cite{GIK} for certain
polynomials on the bidisk $\D^2=\D\times \D$ (here $\D$ is the open unit
disk in the complex plane $\C$).  Define $D(z) = z_1D_1+z_2D_2$ where
the $D_1,D_2$ are $(n+m)\times (n+m)$ matrices given by
\[
D_1 = \begin{pmatrix} I_n & 0 \\ 0 & 0 \end{pmatrix} \quad D_2
= \begin{pmatrix} 0 & 0 \\ 0 & I_m \end{pmatrix}. 
\]

For $n=n_1+n_2$, define
\[
P_{+} = \begin{pmatrix} I_{n_1} & 0 & 0 \\ 0 & 0 & 0 \\ 0 & 0 &
  0 \end{pmatrix} \qquad P_{-} = \begin{pmatrix} 0 & 0 & 0 \\ 0 &
  I_{n_2} & 0 \\ 0 & 0 & 0 \end{pmatrix}
\]
where the blocks correspond to the orthogonal decomposition $\C^{n+m}
= \C^{n_1} \oplus \C^{n_2} \oplus \C^{m}$.  Let $\E = \{z\in \C:
|z|>1\}, \T = \{z\in \C: |z|=1\}$.

\begin{thmb} Suppose $p \in \C[z_1,z_2]$ has bidegree $(n,m)$, no
  zeros in $(\T\times \D) \cup (\T \times \E)$, and no factors depending
  on $z_1$ alone.  Let $n_2$ be the number of zeros of $p(z_1,0)$ in
  $\D$.  Then, there exists an $(n+m)\times (n+m)$ unitary $U$ and a
  constant $c\in \C$ such that
\[
p(z_1,z_2) = c \det ((z_1P_{-}+P_{+} + D_2) - U(P_{-} + z_1P_{+} + z_2
D_2)).
\]
\end{thmb}

This is referred to as a determinantal representation for
``generalized distinguished varieties'' in \cite{GIK} since it
generalizes a determinantal representation for the ``distinguished
varieties'' of Agler and McCarthy \cite{AM} which correspond to the
case $n_2=0$. Polynomials defining distinguished varieties are
essentially a Cayley transform of real stable polynomials and
distinguished varieties have their own motivation in terms of operator
theory as shown in \cite{AM}.  Theorem B is based on first proving a
sums of squares decomposition for polynomials $p \in \C[z_1,z_2]$ with
no zeros in $\T\times \D$ (``a face of the bidisk'') and no factors in
common with $\tilde{p}(z) = z_1^n z_2^m\overline{p(1/\bar{z}_1,
  1/\bar{z_2})}$.  Namely,
\[
|p(z)|^2-|\tilde{p}(z)|^2 = (1-|z_1|^2)(|E_1(z)|^2 - |E_2(z)|^2)  +
(1-|z_2|^2)|F(z)|^2
\]
where $E_1 \in \C^{n_1}[z], E_2 \in \C^{n_2}[z], F \in \C^m[z]$,
$n=n_1+n_2$ where $n_2$ is the number of zeros of $p(z_1,0)$ in $\D$.
This formula generalizes a sums of squares formula of Cole and Wermer
\cite{CW} related to And\^{o}'s inequality from operator theory (see
also \cite{GW} and \cite{pnoz}).  It
would be interesting to characterize which polynomials possess such a
sums of squares formula where $|F(z)|^2$ is also given by a difference
of squares $|F_1(z)|^2-|F_2(z)|^2$, and---going further---it would be
interesting to see what sort of determinantal representation for real
homogeneous polynomials comes out of the corresponding development
from Theorem B to Theorem \ref{mainthm} presented here.

Beyond trivariate polynomials, there are many results on the existence
or non-existence of determinantal representations.  See \cite{vV12},
\cite{KV}, \cite{NT}, \cite{NPT}, \cite{Branden}, \cite{KPV} for
recent results and convenient summaries of the state of the art.
Vinnikov \cite{vV12} conjectures that hyperbolic polynomials always
divide a hyperbolic polynomial which has a determinantal
representation but with additional requirements placed on the set
where the determinantal polynomial is positive.  Our next theorem
offers a step in the right direction for this conjecture albeit in a
special situation.  A polynomial $p$ is \emph{affine} with respect to
a variable $x_j$ if it has degree one in that variable.

\begin{theorem} \label{fourvar} Let $p\in \R[x_0,x_1,x_2,x_3]$ be
  hyperbolic of degree $d$ with respect to the cone
  $\{(0,v_1,v_2,v_3): v_1,v_2,v_3 >0\}$.  Assume $p$ is affine in
  $x_2$ and $x_3$ and of degree $n$ in $x_1$.  Then, there exists
  $k\leq 2n+4$ and $k \times k$ self-adjoint matrices $A_0,A_1,A_2,
  A_3$ such that $p$ divides
\[
\det(\sum_{j=0}^{3} x_j A_j),
\]
$A_1,A_2,A_3$ are positive semi-definite and $A_1+A_2+A_3= I$.
\end{theorem}

See Section \ref{sec:thmc}. Theorem \ref{fourvar} seems to be one of
the few higher dimensional situations where one gets a determinantal
representation from simple hypotheses.  The recent article of Kummer
\cite{Kummer} proves the interesting result that a hyperbolic
polynomial in $n$ variables with no real singularities divides a
determinantal polynomial.  This article also obtains bounds on the
sizes of the matrices involved under the assumption that some power of
the polynomial has a determinantal representation.  Theorem
\ref{fourvar} requires no assumptions of smoothness and obtains
general bounds on the sizes of the matrices involved, but Kummer's
result has the advantage that it works in $n$ variables and does not
assume degree restrictions.

The key tool for Theorem \ref{fourvar} is the following sums of
squares decomposition from Bickel-Knese \cite{BK}.

\begin{thmc}[Theorem 1.12 of \cite{BK}] Let $p\in \C[z_1,z_2,z_3]$
  have multi-degree $(n,1,1)$ and no zeros on $\overline{\D}^3$.
  Then, there exist column-vector valued polynomials $E_1 \in
  \C^n[z_1,z_2,z_3]$, $E_2,E_3\in \C^2[z_1,z_2,z_3]$ such that for
  $z=(z_1,z_2,z_3), w= (w_1,w_2,w_3)$
\[
p(z)\overline{p(w)} - \tilde{p}(z) \overline{\tilde{p}(w)} =
\sum_{j=1}^{3} (1-z_j\bar{w}_j)E_j(w)^* E_j(z)
\]
where $\tilde{p}(z) = z_1^nz_2 z_3
\overline{p(1/\bar{z}_1,1/\bar{z}_2,1/\bar{z}_3)}$.
\end{thmc}

\section{Proof of Theorem \ref{mainthm} from Theorem
  B} \label{sec:btomain} 

Let $\C_{+} = \{z \in \C: \Im z >0 \}, \C_{-} = \{z \in \C: \Im z < 0\}$.

Assume $P \in \R[x_0,x_1,x_2]$ is homogeneous of degree $d$ and for
every $x\in \R^3$
\[
t \mapsto P(x-te_2)
\]
is either identically zero or only has real zeros.  We will assume $P$
has no factors depending only on $x_0,x_1$ which can easily be
incorporated into our final determinantal representation by appending
diagonal blocks to our matrices.

Consider 
\[
q(z_1,z_2) = P(1,z_1,z_2)
\]
which has no zeros in $(\R\times \C_{+}) \cup (\R \times \C_{-})$.  To
see this, take
$z=(a_1,a_2+ib_2) \in (\R\times \C_{+})\cup (\R\times \C_{-})$ with
$q(z)=0$.  Then, $P((1,a_1,a_2) + te_2)$ has the imaginary root
$t=ib_2$, which would imply $t\mapsto P(1,a_1, a_2+t)$ is identically
zero.  This means $x_1-a_1x_0$ divides $P$ which we have ruled out.

Now, define
\[
p(z_1,z_2) = q\left(i\frac{1+z_1}{1-z_1}, i \frac{1+z_2}{1-z_2}\right)
\left(\frac{1-z_1}{2i}\right)^n
\left(\frac{1-z_2}{2i}\right)^m
\]
where $q$ has degree $n$ in $x_1$ and degree $m$ in $x_2$.  Setting
$x_0=1$ in $P(x_0,x_1,x_2)$ cannot lower the degree in $x_1$ or $x_2$,
so $n=\deg_1 P$, $m=\deg_2 P$.  Recall that
\[
z \mapsto i \frac{1+z}{1-z}
\]
is a conformal map of the unit disk onto the upper half plane sending
$\T$ to $\R\cup\{\infty\}$ where $1 \mapsto \infty$.  Thus, $p$ has no zeros in
$(\T\setminus\{1\})\times \D$ as well as $(\T\setminus \{1\}) \times
\E$ where $\E = \{z \in \C: |z|>1\}$.  We cannot have $p(1,z_2) = 0$
unless $p(z_1,z_2)$ has $z_1-1$ as a factor.  This follows by
Hurwitz's theorem since the polynomials $z_2\mapsto p(z_1,z_2)$ will
have no zeros in $\C\setminus\T$ for $z_1 \in \T$ with $z_1 \to 1$, and
then $p(1,z_2)$ will either have the same property or will be
identically zero.  However such factors cannot exist since they imply
$q$ has degree less than $n$ in $x_1$.  In any case, we can safely
divide out factors of $p$ that depend only on $z_1$ since these can
easily be incorporated into our final determinantal representation.
Having done this, $p$ satisfies the hypotheses of Theorem B and we may
write
\[
p(z_1,z_2) = c \det ((z_1P_{-}+P_{+} + D_2) - U(P_{-} + z_1P_{+} + z_2
D_2))
\]
for a unitary $U$. Notice $n_2$ is the number of roots of $z_1\mapsto
p(z_1,0)$ in $\D$ which is the same as the number of roots of $z_1
\mapsto q(z_1, i)=P(1,z_1,i)$ in $\C_{+}$.

We convert back to $q$ via $z \mapsto \frac{z-i}{z+i}$.
So,
\begin{align}
  q(z_1,z_2) &=  p\left(\frac{z_1-i}{z_1+i},
    \frac{z_2-i}{z_2+i}\right) (z_1+i)^n (z_2+i)^m \nonumber \\
&= p\left(\frac{z_1-i}{z_1+i},
    \frac{z_2-i}{z_2+i}\right) \det( (z_1+i) D_1 + (z_2+i)D_2) \nonumber\\
&= \begin{aligned} c\det( (z_1-i)P_{-} &+ (z_1+i)P_{+} + (z_2+i)D_2\\
 &- U((z_1+i)P_{-} +
(z_1-i)P_{+} +(z_2-i)D_2))\end{aligned} \nonumber \\
&= c\det( (I-U) D(z) - i(I+U)(P_{-} - P_{+} - D_2)) \nonumber \\
&= \pm c \det( (I-U)(-z_1P_{-} + z_1 P_{+} +z_2 D_2) + i(I+U)). \label{detcomps} 
\end{align}
The last line comes from multiplying on the right by $\det(-P_{-}
+P_{+} + D_2)$.  Letting $M(z) = -z_1P_{-} + z_1 P_{+} +z_2 D_2$, we
now form the spectral decomposition $U= V \begin{pmatrix} u & 0 \\ 0 &
  I \end{pmatrix} V^*$; $V$ is a unitary, $u$ is a $k\times k$
diagonal unitary with no $1$'s on the diagonal, and $k$ is the rank of
$U-I$.  Factoring $V$ and $V^{*}$ out from the left and right of
\eqref{detcomps} leaves
\[
\begin{aligned}
q(z) =&\pm c \det \left( \begin{pmatrix} I-u & 0 \\ 0 & 0 \end{pmatrix} V^*
  M(z) V + i \begin{pmatrix} I+u & 0 \\ 0 & 2I \end{pmatrix} \right)
\\
&=\pm c\det(I-u) \det \left( \begin{pmatrix} I & 0 \\ 0 & 0 \end{pmatrix} V^*
  M(z) V + \begin{pmatrix} a & 0 \\ 0 & 2iI \end{pmatrix} \right)
\\ 
&= \pm c\det(I-u) \det \begin{pmatrix} (V^*M(z)V)_{kk} + a & * \\ 0 &
  2iI \end{pmatrix} \\
&= C \det((V^* M(z) V)_{kk}+a)
\end{aligned}
\]
where $a = i(I+u)(I-u)^{-1}$ is a diagonal matrix with real entries,
$(V^*M(z)V)_{kk}$ is the upper $k\times k$ block of $V^*M(z)M$, and
$C$ is a constant.  Now, $V^*M(z)V = -z_1V^*P_{-}V + z_1 V^* P_{+}V
+z_2 V^*D_2V$ and if we set $A_0 =a$, $A_1 =
(-V^*P_{-}V + V^*P_{+}V)_{kk}$, and $A_2 = (V^*D_2V)_{kk}$ we have 
a determinantal representation for $q$: 
\[
q(z) = C \det( A_0+z_1 A_1 + z_2 A_2).
\]
Notice $A_0,A_1,A_2$ are evidently self-adjoint with $A_2$ positive
semi-definite, and since $\deg q = d$ we have $d \leq k$.  Once we
show $k=d$, we can homogenize to get the determinantal representation
for $P$.  It helps to first establish some of the additional details
listed in Theorem \ref{mainthm}.

It is a general fact that for matrices $A,B$, the degree of
$\det(tA+B)$ is at most $\rank A$ (we leave this as an exercise).  So,
$\deg_j q \leq \rank A_j$ for $j=1,2$.  On the other hand, by
construction $\rank A_1 \leq \rank (-P_{-}+P_{+}) = \deg_1 q$ and
$\rank A_2 \leq \rank D_2 = \deg_2 q$, yielding $\deg_j q = \rank A_j$
for $j=1,2$.  Next, setting $B_{\pm} = (V^*P_{\pm} V)_{kk}$ we have
$A_1 = B_{+}-B_{-}$.  Since $\rank A_1 = n_1+n_2$ and $\rank B_{+}
\leq n_1$ and $\rank B_{-} \leq n_2$, we must have equality in both
inequalities.  This also shows the ranges of $B_{+},B_{-}$ have
trivial intersection by considering dimensions.
Since $P_{+}+P_{-} + D_2 = I$, we must have
$B_{+}+B_{-} + A_2 = I$. 

In order to show $k=d$, it suffices to show $Q(t):=t A_1+A_2$ is
non-singular for some $t$.  For then, there would be a $t_0$ such that
\[
t \mapsto q(t(t_0,1))
\]
has degree $k$, and since $q$ has degree $d$, we would have $k\leq d$
and thus $k=d$.

Note $Q(t) = I +(t - 1)B_{+} - (t+1)B_{-}$.  By the spectral theorem 
\[
B_{+} = \sum_{j=1}^{n_1} \nu_j v_j v_j^* \qquad
B_{-} = \sum_{j=1}^{n_2} \mu_j w_j w_j^*
\]
where $V=\{v_1,\dots, v_{n_1}\}, W=\{w_1,\dots, w_{n_2}\}$ form
orthonormal sets of eigenvectors corresponding to the positive
eigenvalues $\{\nu_1,\dots, \nu_{n_1}\}$,  $\{\mu_1,\dots, \mu_{n_2}\}$
of $B_{+}$, $B_{-}$ respectively.  Let $Y=\{y_1,\dots y_{k-n}\}$ be an
orthonormal basis for the complement of $V$ and $W$.  Then,
$\mathcal{B} = V\cup W\cup Y$ is a basis for $\C^k$.  Let
$\mathcal{C}$ be a basis dual to $\mathcal{B}$.  (Two bases
$\{b_1,\dots b_N\}$, $\{c_1,\dots, c_N\}$ are \emph{dual} if $b_j^*c_k
= \delta_{jk}$.)  The matrix for $Q(t)$ obtained by using
$\mathcal{C}$ as a basis for the domain and $\mathcal{B}$ for the
range is of the form
\[
\begin{pmatrix} I +(t-1)d_{+} & 0 & 0 \\ 0 & I-(t+1)d_{-}& 0 \\ 0 & 0 & I 
\end{pmatrix}
\]
for diagonal matrices $d_{+},d_{-}$ containing the eigenvalues
$\nu_1,\dots,\nu_{n_1}$,$\mu_1,\dots, \mu_{n_2}$ on the diagonal.  The
determinant of this vanishes for only finitely many $t$ and so
$Q(t_0)$ is certainly non-singular for some $t_0$.  Thus, $k=d$ and we
homogenize $q$ at degree $d$ to see that
\[ 
P(x) = C\det(x_0A_0+ x_1A_1 + x_2 A_2).
\]
This concludes the proof of Theorem \ref{mainthm}.

\section{Example} \label{sec:ex} 
Renegar \cite{Renegar} has an example of a polynomial that is
semi-hyperbolic but not hyperbolic in any direction  (see Section
2 of that paper); however we have constructed an example that is more
illustrative for our purposes.

Let
\[
p(x_0,x_1,x_2) = 2x_0^2 x_1-(x_0^2+3x_1^2)x_2.
\]
Then, $t \mapsto p(x-te_2)$ clearly has only real roots for $x \in
\R^3$ since this one variable polynomial has degree $1$ and real
coefficients.  Let
\[
A_0 = \frac{i}{3}\begin{pmatrix} 0 & -3 & -\sqrt{3} \\
3 & 0 & \sqrt{3} \\ \sqrt{3} & -\sqrt{3} & 0 \end{pmatrix} 
\quad
A_1 =\begin{pmatrix} 0 & 0 & 0 \\ 0 & 1 & 0 \\ 0 & 0 & -1 \end{pmatrix} 
\quad
A_2 = \begin{pmatrix} 1 & 0 & 0 \\ 0 & 0 & 0 \\ 0 & 0 & 0\end{pmatrix}.
\]
We see that
\[
p(x) = 3\det(x_0A_0 +x_1A_1+x_2A_2).
\]
As remarked in the introduction we can lift to
\[
P(x_0,x_1,y_1,x_2) = 3x_1y_1 x_2 - (x_2+ x_1 + 3y_1)x_0^2
\]
which is hyperbolic in the direction $(0,1,1,1)$ and
$P(x_0,x_1,-x_1,x_2) = p(x_0,x_1,x_2)$.  We now explain why $p$ is not
hyperbolic in any direction.

We first show that $\{x:p(x)\ne 0\}$ consists of the two connected
components $\mcP_{+} = \{x:p(x)>0\}$, $\mcP_{-}=\{x:p(x)<0\}$.  I
thank the referee for the following simplified explanation.  The
hypersurface $\{x : p(x) = 0\}$ is the graph of the continuous
function $(x_0,x_1) \mapsto \frac{2x_0^2x_1}{x_0^2 + 3x_1^2}$.  Thus,
$\{x:p(x) \ne 0\}$ is divided into exactly two components: the part
above the graph and the part below.


Next, neither component $\mcP_{+}, \mcP_{-}$ is convex.   
For instance, $(-1,0,-1)$, $(1,0,-1)\in P_{+}$ but $(0,0,-1) \notin \mcP_{+}$.
One can similarly show $\mcP_{-}$ is not convex.  This implies that
$p$ is not hyperbolic in any direction since it is a fundamental result
of G\r{a}rding that if $p$ is hyperbolic in some direction $e$, then the 
connected component of $\{x:p(x)\ne 0\}$ containing $e$ is convex.  

This brings up a potential paradox.  Since $p$ is a limit of
hyperbolic polynomials $p_{\epsilon}$ as in equation \eqref{pepsilon}, 
how is it possible that the connected components of
$\{x:p(x)\ne 0 \}$ are non-convex in the above example?  \emph{An}
answer is that a convex component of $\{x:p_\epsilon(x)\ne 0\}$ could
shrink to an isolated point (in projective space) as $\epsilon
\searrow 0$.  This is something we have seen graphically using the
above example.

Finally, in connection with our discussion after Theorem \ref{mainthm}
regarding the signatures of $A_0,A_1$, let us point out that
\[
p(1,t,i) = 2t-(1+3t^2)i
\]
has one zero in $\C_{+}$, which agrees with the number of negative
eigenvalues of $A_1$.  On the other hand, $p(0,t,i) = -3t^2i$ has no
zeros in $\C_{+}$.  Notice also that $p(t,1,i) = (2-i)t^2-3i$ has one
zero in $\C_{+}$.  This matches the number of negative eigenvalues of
$A_0$, which is what one would like to have more generally in order
for Theorem \ref{mainthm} to have a more symmetric statement.

\section{Proof of Corollary \ref{cor}} \label{sec:cor} Notice that
$t \mapsto p(x-te_2)$ is either identically zero or only has real
roots by Hurwitz's theorem since this polynomial can be obtained as
the limit as $a \searrow 0$ of
\[
t \mapsto p(x-t(ae_1+e_2)).
\]
Any factors depending only on $x_0,x_1$ can easily be dealt with
separately so we may assume there are no such factors.
So, $p$ satisfies the hypotheses of Theorem \ref{mainthm}.
Also, $t\mapsto p(1,t,i)$ can have no zeros in the upper half plane 
for if it had such a zero $z = x+iy$ where $y>0$, then
\[
t \mapsto p((1,x,0)+t(0,y,1))
\]
would have the non-real zero $t=i$ contradicting hyperbolicity in the
direction $(0,y,1)$.  This shows that $\rank B_{-}=0$ in Theorem \ref{mainthm}
and therefore $A_1$ is positive semi-definite as desired.

\section{Proof of Theorem \ref{fourvar} from Theorem C}\label{sec:thmc}
We largely follow the scheme of \cite{gK}.
Let $P\in \R[x_0,x_1,x_2,x_3]$ be homogeneous of degree $d$ of degree
$1$ in $x_2,x_3$ and of degree $n$ in $x_1$.  Assume $P$ is hyperbolic
with respect to the cone $\{(0,v_1,v_2,v_3): v_1,v_2,v_3>0\}$.
Then, for $x=(x_1,x_2,x_3)$
\[
q(x) = P(1,x)
\]
has no zeros in $\C_{+}^3 \cup \C_{-}^3$ and
$\overline{q(\bar{x})}=q(x)$.  Switching to the tridisk, we see that
\[
f(z) = q\left(i\frac{1+z_1}{1-z_1}, i \frac{1+z_2}{1-z_2},
  i\frac{1+z_3}{1-z_3}\right)\left(\frac{1-z_1}{2i}\right)^n\left(\frac{1-z_2}{2i}\right)
\left(\frac{1-z_3}{2i}\right)
\] 
has no zeros in $\D^3\cup \E^3$.  Note that we may as well assume $f$
is irreducible since otherwise $f$ will have a factor depending on one
or two variables alone, in which case there is no issue with having a
determinantal representation.

Let $1/\bar{z} = (1/\bar{z}_1,1/\bar{z}_2,1/\bar{z}_3)$ for $z \in
\C^3$ and define
\[
\begin{aligned}
\tilde{f}(z) &= z_1^n z_2 z_3 \overline{f(1/\bar{z})} \\
  \widetilde{\frac{\partial
    f}{\partial z_j}} &= \frac{z_1^{n}z_2z_3}{z_j} \overline{\frac{\partial
    f}{\partial z_j}(1/\bar{z})} \quad \text{ for } j=1,2,3.
\end{aligned}
\]

Since $q$ has real coefficients one can show that $\tilde{f}=f$ and
\[
\begin{aligned}
nf &= z_1 \frac{\partial f}{\partial z_1} +  \widetilde{\frac{\partial
    f}{\partial z_1}} \\
f &= z_j \frac{\partial f}{\partial z_j} + \widetilde{\frac{\partial
    f}{\partial z_j}} \text{ for } j=2,3
\end{aligned}
\]
after some simple computations. Thus, $(n+2)f = p + \tilde{p}$ where
\[
p(z) = \sum_{j=1}^{3}\widetilde{\frac{\partial f}{\partial z_j}} \qquad 
\tilde{p}(z) = \sum_{j=1}^{3} z_j \frac{\partial f}{\partial z_j}.
\]
Let $f_t(z) = f(tz)$ for $0<t<1$.  Then, $f_t$ has no zeros in
$\overline{\D}^3$ and if we set $\tilde{f_t}(z) = t^{n+2}f(z/t)$, then
$|f_t|=|\tilde{f_t}|$ on $\T^3$ (since $\tilde{f}=f$) and so
$\tilde{f_t}/f_t$ is analytic and bounded by $1$ in modulus for $z\in
\D^3$ by the maximum principle.  Now, for $z \in \cD^3$
\[
\begin{aligned}
0 &\leq \lim_{t\nearrow 1} \frac{|f(tz)|^2 - |t^{n+2}f(z/t)|^2}{1-t^2}(n+2)
\\
&= (n+2)^2|f(z)|^2 -2 \text{Re}(\tilde{p}(z)(n+2)\overline{f(z)})\\
& = |p(z)|^2
- |\tilde{p}(z)|^2 \text{ since } (n+2)f=p+\tilde{p}
\end{aligned}
\]
with some computations omitted (see \cite{gK} for more details).  This
shows that if $p$ vanishes in $\D^3$, then so does $\tilde{p}$ and so
does $f$ which by assumption does not happen.  Hence, $p$ has no zeros
in $\D^3$.

Note that if $p$ and $\tilde{p}$ had a common factor then this would
be a factor of $f$ which we have already ruled out; we point out that
$p$ and $\tilde{p}$ cannot be multiples of one another since
$\tilde{p}$ vanishes at the origin.  The \emph{conclusion} of Theorem
C holds for such a $p$ but since we have only stated it for
polynomials with no zeros on $\cD^3$ (as opposed to $\D^3$) we must
explain how to address the case at hand.  The main point is that for
$0<t<1$, $p_t(z) = p(tz)$ will satisfy the hypotheses of Theorem C and
therefore there exist vector polynomials $E_1^t,E_2^t, E_3^t$
corresponding to $p_t$ as in Theorem C.  Then,
\[
\sup_{\T^3} |p(z)|^2 \geq (1-|z_j|^2)|E_j^t(z)|^2
\]
shows the vector polynomials $E^t_j$ are locally bounded in $\D^3$ and
hence we can choose subsequences of $t\nearrow 1$ such that that
$E_1^t\in \C^{2n}[z],E_2^t, E_3^t \in \C^2[z]$ converge to vector
polynomials $E_1\in \C^{2n}[z],E_2,E_3\in \C^2[z]$ and hence we will
get a sums of squares decomposition as in Theorem C.  Note the
polynomials in $E_1,E_2,E_3$ necessarily have degree at most
$(n-1,1,1), (n,0,1),(n,1,0)$ (this is proven in \cite{gKrifitsac} for
instance) and they will be non-trivial since $p$ and $\tilde{p}$ have
no factors in common.  On the zero set $Z_f$ of $f$, $p=-\tilde{p}$
and therefore
\begin{equation} \label{sosonz}
0 = \sum_{j=1}^{3} (1-z_j\bar{w}_j)E_j(w)^* E_j(z)
\end{equation}
for $z,w \in Z_f$.  This equation ensures that the map 
\begin{equation} \label{lurkingisom}
\begin{pmatrix} z_1 E_1(z) \\ z_2 E_2(z) \\ z_3 E_3(z) \end{pmatrix}
\mapsto \begin{pmatrix}  E_1(z) \\  E_2(z) \\  E_3(z) \end{pmatrix}
\end{equation}
defined initially for vectors of the above form with $z\in Z_f$,
extends linearly to a well-defined $(2n+4)\times (2n+4)$ unitary $U$.
(Some details: If a combination of vectors from the left side of
\eqref{lurkingisom} sums to zero, \eqref{sosonz} shows the
corresponding combination on the right sums to zero. So, we get a
well-defined linear map from the span of the left side of
\eqref{lurkingisom} to the span of the right side.  Now,
\eqref{sosonz} shows this map is an isometry.  Since we are in finite
dimensions it can be extended to a unitary.)

Note that $E_1,E_2,E_3$ cannot vanish identically in $Z_f$ without
vanishing in all of $\C^3$ since the degrees are lower and $f$ is
irreducible.  Let $P_j$ for $j=1,2,3$ be the projection onto the
$j$-th component in the orthogonal decomposition of $\C^{2n+4} =
\C^{2n}\oplus \C^{2} \oplus \C^2$ and let $M(z) = \sum_{j=1}^3 z_j P_j$.
By \eqref{lurkingisom}, for $z \in Z_f$
\[
(I-U M(z)) \begin{pmatrix} E_1(z) \\ E_2(z) \\ E_3(z) \end{pmatrix} =
0
\]
and therefore $\det(I-UM(z)) = 0$ for $z \in Z_f \setminus \{z:
E_1,E_2,E_3=0\}$.  Basic results in algebraic geometry (such as in
Chapter 4, Section 4 of \cite{IVA}) can be used to establish that this
implies $\det(I-UM(z))$ vanishes for $z \in Z_f$ (i.e. $Z_f \setminus
\{z: E_1,E_2,E_3=0\}$ is Zariski dense in $Z_f$) since $f$ is
irreducible and none of $E_1,E_2,E_3$ vanish identically on $Z_f$.

Therefore $f$ divides $\det(I-UM(z))$.  Write
\[
f(z) g(z) = \det( I - U M(z))
\]
for some polynomial $g$ of degree at most $(n,1,1)$. As with
Section \ref{sec:btomain}, we convert back to $q$.  There is some
repetition in what follows but since the situations are slightly
different we include the details. Now,
\begin{align}
q(z) r(z) &= \det( (\sum_{j=1}^3 (z_j+i)P_j) - U (\sum_{j=1}^3 (z_j-i)P_j) ) 
\nonumber\\
&= \det( (I-U)M(z)+i(I+U)) \label{qrdet}
\end{align}
for 
\[
r(z) = (z_1+i)^n (z_2+i)(z_3+i) g\left(\frac{z_1-i}{z_1+i},
    \frac{z_2-i}{z_2+i}, \frac{z_3-i}{z_3+i}\right).
\]
Let $U= V \begin{pmatrix} u & 0 \\ 0 & I \end{pmatrix} V^*$ be the
spectral decomposition of $U$ where $u$ is $k\times k$ diagonal with
unimodular entries, none of which equals $1$. Here $k$ is the rank of
$I-U$.  As in Section \ref{sec:btomain}, the determinant \eqref{qrdet}
can be converted to
\begin{equation} \label{convertedqr}
q(z) r(z) = (\text{const}) \det ( (V^* M(z) V)_{kk} + a)
\end{equation}
where again $(V^*M(z)V)_{kk}$ refers to taking the upper $k\times k$
block of the given matrix, and $a = i(I+u)(I-u)^{-1}$.  Finally, if we
homogenize \eqref{convertedqr} at degree $k$---note this is at most
$2n+4$---then
\[
P(x) R(x) = (\text{const}) \det( x_0 a + \sum_{j=1}^{3} x_j A_j)
\]
with $A_j = (V^*P_jV)_{kk}$ and $A_1+A_2+A_3 = (V^*IV)_{kk} = I$, and
where $R(x) = x_0^{k-d} r((1/x_0)(x_1,x_2,x_3))$.  This concludes the
proof.

\section{Concluding questions and remarks} \label{conclusion}
We think it is worthwhile to discuss or rehash some of the motivations
and lingering questions of this paper in more detail.

Semi-hyperbolic polynomials have perhaps been overlooked because they
lack one of the key features of hyperbolic polynomials.  Specifically,
if $p$ is hyperbolic in the direction $e$, then the connected
component of $\{x:p(x)\ne 0\}$ containing $e$ is convex (see
\cite{Renegar}).  No such result holds for semi-hyperbolic polynomials
(see Renegar \cite{Renegar} Section 2 or Section \ref{sec:ex} above).
This convexity property ties hyperbolic polynomials to optimization
and is ``the cornerstone of hyperbolic programming'' \cite{Renegar}.
This begs the question, why study \emph{semi}-hyperbolic polynomials
which may lack this property?

First, we think it is a good general principle in mathematics to
understand the degenerate versions of objects of interest.  Notice
that the (local uniform) limit of a sequence of homogeneous
polynomials of degree $d$ which are semi-hyperbolic with respect to a
specific direction $e$ is either semi-hyperbolic or identically zero.
This follows from Hurwitz's theorem applied to each polynomial
$t\mapsto P(x-te)$.  Hyperbolic polynomials do not share this
property.  Somewhat related is the following question mentioned in the
introduction.

\begin{question} Every trivariate semi-hyperbolic polynomial is a
  limit of hyperbolic polynomials.  Is this true more variables?
\end{question}

Second, our main theorem, Theorem \ref{mainthm}, shows that trivariate
semi-hyperbolic polynomials possess determinantal representations just
as in the hyperbolic case.  We think this in itself provides good
justification for the study of semi-hyperbolic polynomials.  This is a
good point to formally state a question from the introduction.  

\begin{question} In Theorem \ref{mainthm}, can the signature and rank
  of $A_0$ be determined directly from properties of $p$?
\end{question}

Our own personal motivations for studying semi-hyperbolic polynomials
came from the natural connection we presented above between
semi-hyperbolic polynomials and two variable polynomials with no zeros
on $\T\times \D$.  These latter polynomials appeared in \cite{GIK}
essentially because of the realization that some of the theory of
polynomials with no zeros in $\D^2$ could be pushed further to the
situation of no zeros in $\T\times \D$.  It was realized later that
this initially unnatural condition is closely related to hyperbolicity
and indeed is essentially equivalent to semi-hyperbolicity.  

Finally, we wish to rehash our larger goal of the paper of connecting
sums of squares formulas to determinantal representations. This
description will be somewhat imprecise. The approach of this paper
shows that if $p(z_1,z_2,\dots,z_n)$ has no zeros in $\D^n$ and
possesses a hermitian sums of squares formula
\begin{equation} \label{agler}
|p|^2 - |\tilde{p}|^2 = \sum_{j=1}^{n} (1-|z_j|^2)SOS_j
\end{equation}
(here each $SOS_j$ term is a sum of squared moduli of polynomials)
then $p+\tilde{p}$ divides a unitary determinantal polynomial
\[
\det(I-U D(z)).
\]
Here $U$ is a unitary matrix and $D(z)$ is a diagonal matrix with
coordinate functions on the diagonal.  One can then convert
$p+\tilde{p}$ via Cayley transform and homogenization to a hyperbolic
polynomial (hyperbolic with respect to all vectors with positive
entries) and through some linear algebra get a self-adjoint
determinantal polynomial.  One can reverse engineer some of this: take
a hyperbolic polynomial $P$ (again hyperbolic with respect to vectors
with positive entries) convert to a polynomial $q$ satisfying
$q=\tilde{q}$.  If $q$ can be written as $p+\tilde{p}$ where $p$
satisfies \eqref{agler} then $P$ divides a determinantal
representation.  If $n>2$, not every $p\in \C[z_1,\dots,n]$ with no
zeros in $\D^n$ satisfies an equation of the form \eqref{agler}; such
polynomials are called \emph{Agler denominators}.  A polynomial is an
Agler denominator if and only if $\tilde{p}/p$ satisfies a
multivariable von Neumann inequality (see \cite{gKrifitsac}).  

We have also presented a modification to
hyperbolicity/semi-hyperbolicity with respect to a specific direction.
In $n$ variables this would entail, after various conversions, to
understanding polynomials satisfying \eqref{agler} where
$SOS_1,\dots SOS_{n-1}$ are replaced with differences of squares.
With the exception of our work in \cite{GIK}, this is relatively
uncharted territory.

\section*{Acknowledgments}

I would like to sincerely thank my co-authors J. Geronimo and P. Iliev
\cite{GIK}, my co-author K. Bickel \cite{BK}, and the anonymous
referee.  I would like to also thank Victor Vinnikov for introducing
me to hyperbolic polynomials.

\end{document}